\documentclass[11pt]{amsart}

\usepackage{amssymb, graphicx}
\usepackage{graphicx}

\usepackage[all]{xy}

\usepackage{amsfonts}
\usepackage{amsmath}
\usepackage{latexsym}
\usepackage{amscd}
\usepackage{verbatim}

\allowdisplaybreaks[1]

\newtheorem{theorem}{Theorem}[section]
\newtheorem{proposition}[theorem]{Proposition}

\theoremstyle{remark}

\begin{document}
\title[On the Theorem of Amitsur--Levitzki]{On the Theorem of Amitsur--Levitzki}

\author{Claudio Procesi}

\address{C. Procesi, Dipartimento di Matematica, Sapienza
 Universit\` a di Roma,  Italy}  \email{procesi@mat.uniroma1.it}

\begin{abstract} We present a proof of the Amitsur--Levitzki theorem which is  a basis for a general theory of equivariant skew--symmetric maps on matrices.\medskip

\hskip4cm{\em to the memory of  S.A. Amitsur who taught me PI theory}\end{abstract}

\thanks{2010 {\em Math. Subj. Class.} Primary 16R60. Secondary 16R10, 16R30. }
\maketitle
\section{Introduction}
The Theorem of Amitsur--Levitzki is a cornerstone of the Theory of polynomial identities.  It states that the algebra of $n\times n$ matrices over any commutative ring $A$  satisfies the standard polynomial $S_{2n}$.   For any integer $h$   the standard polynomial $S_{h}$ is the element of the free algebra  in the variables $x_1,\ldots,x_h$ given by the formula
$$S_h(x_1,\ldots,x_h):=\sum_{\sigma\in \mathfrak S_h}\epsilon_\sigma x_{\sigma(1)}\ldots x_{\sigma(h)} $$ where $\mathfrak S_h$ is the symmetric group of permutations on $\{1,\ldots,h\}$ and $\epsilon_\sigma$  denotes the sign of the permutation $\sigma$.\smallskip

This Theorem has received several proofs, after the original proof \cite{AmL}, which is a direct verification. A  proof by Swan using graphs, \cite{Sw}; one by Kostant  relating it to Lie algebra cohomology, \cite{kost}; by   Formanek as consequence of the Cayley--Hamilton Theorem (see \cite{Fo1} and  \cite{Fo2}) and finally a very simple by Rosset using Grassman variables \cite{Ros}.  So why give another proof? 

 The proof I present shows that the Amitsur--Levitzki Theorem is the Cayley--Hamilton identity for the {\em generic Grassman matrix}. Techinally it  is very similar   to Rosset's proof but while in his proof the Grassman variables  are auxiliary, in the present proof these variables are intrinsically embedded in the problem, this is important for applications.  
 
 In this formulation the Theorem is the first step for the general Theory  of alternating equivariant maps (with respect to conjugation)  from matrices to matrices.  
This theory can be considered as the Grassman analogue of the theory of generic matrices with trace and it  is fully explained in a joint paper with Matej Bresar and Spela Spenko  \cite{bps}  but I think it is worth presenting the result on the Amitsur--Levitzki Theorem independently.

I will state at the end of the paper the general Theorem \ref{alg} from \cite{bps} which shows that this algebra is very explicit.  

As usual it is enough to prove the Theorem for matrices over a field $\mathbb F$ of characteristic 0, i.e. $\mathbb Q$ so we will assume we are in this setting from now on.

\section{Antisymmetry} \label{assec4}

\subsection{Antisymmetry}
By the {\em antisymmetrizer} we mean the operator that sends a multilinear expression $f(x_1,\ldots,x_h)$ into the antisymmetric expression $ \sum_{\sigma\in\mathfrak  S_h} \epsilon_\sigma f(x_{\sigma(1)},\ldots,x_{\sigma(h)})$; so that. if $\phi_1,\ldots,\phi_h$ are linear forms on $V$,  antisymmetrizing $\phi_1(v_1)\ldots\phi_h(v_h)$ one has $\phi_1\wedge\ldots\wedge\phi_h$ and, applying the antisymmetrizer to the noncommutative monomial $x_1\cdots x_h$ we get  
the standard polynomial of degree $h$,  $S_h(x_1,\ldots,x_h)=\sum_{\sigma\in\mathfrak  S_{h}}\epsilon_\sigma   x_{\sigma(1)} \dots x_{\sigma(h)}.$
Up to a scalar multiple this is the only multilinear antisymmetric noncommutative polynomial of degree $h$.

Let $A$ be any $\mathbb F$-algebra (not necessarily associative) with basis $e_i$,  and let $V$  be a finite dimensional vector space  over $\mathbb F$. The set of  multilinear  antisymmetric functions from $V^k$ to $A$  is given by functions $G(v_1,\ldots,v_k)=\sum_iG_i(v_1,\ldots,v_k)e_i$  with $G_i(v_1,\ldots,v_k)$ multilinear  antisymmetric functions from $V^k$ to $\mathbb F$,  moreover if $A$ is infinite dimensional only finitely many $G_i(v_1,\ldots,v_k)$ appear for any given $G$.

 In other words  $G_i(v_1,\ldots,v_k)\in \bigwedge^kV^*$ and
 this space can be  identified with $\bigwedge^kV^*\otimes A$.  Using the algebra structure of $A$ we have a wedge product of these functions:
\smallskip

 \noindent for $G\in \bigwedge^hV^*\otimes A,\ H\in \bigwedge^kV^*\otimes A$ we define
$$(G\wedge H)( v_1,\ldots,v_{h+k}):=\frac{1}{h!k!}\sum_{\sigma\in \mathfrak S_{h+k}}\epsilon_\sigma G (v_{\sigma(1)},\ldots,v_{\sigma(h)})H(v_{\sigma(h+1) },\ldots,v_{\sigma(h+k)}) $$
$$=\sum_{\sigma\in \mathfrak S_{h+k}/\mathfrak S_{h}\times \mathfrak S_{k}}\epsilon_\sigma G (v_{\sigma(1)},\ldots,v_{\sigma(h)})H(v_{\sigma(h+1) },\ldots,v_{\sigma(h+k)}).$$
As an example we have $S_a\wedge S_b=S_{a+b} $.

With this multiplication 
the algebra of multilinear  antisymmetric functions from $V $ to $A$ is isomorphic to the tensor product algebra $\bigwedge V^*\otimes A$.  We shall denote by $\wedge$ the product in this algebra.

Assume now that $A$ is an associative algebra and $V\subset A$.  The inclusion map $X:V\to A$    is of course antisymmetric, since the symmetric group on one variable is trivial, hence $X\in \bigwedge V^*\otimes A$. By iterating the definition of wedge product we have the important fact:
\begin{proposition}
As a multilinear function, each power  $X^a:=X^{\wedge a}$  equals the standard polynomial $S_a$  computed in $V$.
\end{proposition}
We apply this to $V=A=M_n(\mathbb F):=M_n$;  the group $G=PGL(n,\mathbb F)$ acts  on this space and hence on functions by conjugation and it is interesting to study the {\em invariant algebra}, i.e. the algebra of $G$--equivariant maps\begin{equation}
\label{an}A_n:=(\bigwedge M_n^*\otimes M_n)^G.
\end{equation} 
This among other topics is  discussed in \cite{bps}.

In the natural basis $e_{ij}$ of matrices and the coordinates $x_{ij}$ the element $X\in A_n$ (cf. \eqref{an}) is  the {\em   generic Grassman matrix}\quad $X=\sum_{h,k}x_{hk}e_{hk} $.

Hence in this language the Amitsur--Levitzki Theorem  is the single identity $X^{2n}=0$.

\begin{proof}[Proof of Amitsur--Levitzki $X^{2n}=0$.]  Since $X$ is an element of degree 1 we have that $X^2$  is in $\bigwedge^{2} M_n^*\otimes M_n\subset M_n(\bigwedge^{even} M_n^*)$.  We have that 
$X^2$ is an $n\times n$ matrix  over a commutative ring,   the even part of the Grassman algebra, hence in order to prove that $X^{2n}=(X^2)^n=0$  we need to show that $tr((X^2)^{i})=tr(X^{2i})=0$  for $i\leq n$.   Now the fact that the trace of an even  standard polynomial vanishes is an easy exercise and it appears also   in Kostant and in Rosset's proof   so I will not reproduce it here. Thus the Theorem is proved.
\end{proof}
As I already mentioned this formulation can be taken as basis of the description of the algebra $A_n$ (Formula \eqref{an}) of equivariant multilinear antisymmetric maps from matrices to matrices. In \cite{bps} we prove among other results
\begin{theorem}
\label{alg}  The algebra $A_n$ is generated by $X$  and the elements $tr(X^{2i-1}), \ i=1,\ldots, n.$  All these elements anti commute.

$A_n$  is a free module with basis $X^i,\ i=0,\ldots, 2n-1$ over the Grassman algebra in the elements $tr(X^{2i-1}), \ i=1,\ldots, n-1$ and we have the two defining identities
$$X^{2n}=0,\quad tr(X^{2n-1})=-\sum_{i=1}^{n-1}X^{2i} \wedge  tr(X^{2(n-i)-1})+nX^{2n-1}.$$
\end{theorem} 
This Theorem is based on the first and second fundamental Theorem for matrix invariants and in turn gives rise to a series of questions on Lie algebras  which will be treated elsewhere.
\bibliographystyle{plain}

\bibliography{bibliografia}

\printindex

\end{document}